\documentclass[english,11pt]{article}
\usepackage{}
\usepackage{amsmath, amssymb, amstext, hyperref}
\usepackage{graphicx, color, framed}
\graphicspath{{./Figures/}}

\newtheorem{theorem}{Theorem}

\newenvironment{proof}[1][Proof]{\textbf{#1.} }{\  \rule{0.5em}{0.5em}}

\DeclareMathOperator*{\argmin}{argmin}

\DeclareMathOperator{\vect}{vec}

\title{\bf Optimization of Structural Similarity in Mathematical Imaging}

\author{D. Otero\footnotemark[1]
\and D. La Torre\footnotemark[2]
\and O. Michailovich\footnotemark[3]
\and E.~R. Vrscay\footnotemark[1] \footnotemark[4]}

\begin{document}
\maketitle

\newcommand\red[1]{{\color{red}#1}}
\newcommand\blue[1]{{\color{blue}#1}}
\renewcommand{\thefootnote}{\fnsymbol{footnote}}

\footnotetext[1]{Department of Applied Mathematics, University of Waterloo, ON, Canada.}
\footnotetext[2]{SKEMA Business School and Universit$\acute{e}$ C$\hat{o}$te d'Azur, Sophia Antipolis Campus, France. Email: davide.latorre@skema.edu}
\footnotetext[3]{Department of Electrical and Computer Engineering, University of Waterloo, ON, Canada.}
\footnotetext[4]{Questions, comments, or corrections to this document may be directed to the following email address: {\tt ervrscay@uwaterloo.ca}.}

\renewcommand{\thefootnote}{\arabic{footnote}}

\begin{abstract}
It is now generally accepted that Euclidean-based metrics may not always adequately represent the subjective judgement of a human observer. As a result, many image processing methodologies have been recently extended to take advantage of alternative visual quality measures, the most prominent of which is the Structural Similarity Index Measure (SSIM). The superiority of the latter over Euclidean-based metrics have been demonstrated in several studies. However, being focused on specific applications, the findings of such studies often lack generality which, if otherwise acknowledged, could have provided a useful guidance for further development of SSIM-based image processing algorithms. Accordingly, instead of focusing on a particular image processing task, in this paper, we introduce a general framework that encompasses a wide range of imaging applications in which the SSIM can be employed as a fidelity measure. Subsequently, we show how the framework can be used to cast some standard as well as original imaging tasks into optimization problems, followed by a discussion of a number of novel numerical strategies for their solution.
\end{abstract}


\pagestyle{myheadings}
\thispagestyle{plain}
\markboth{D. Otero, E.~R. Vrscay, O. Michailovich and D. La Torre}{Optimization of Structural Similarity in Mathematical Imaging}

\section{Introduction}
Image denoising, deblurring, and inpainting are only a few examples of standard image processing tasks which are traditionally solved through numerical optimization. In most cases, the objective function associated with such problems is expressed as the sum of a {\em data fidelity term} $f$ and a {\em regularization term} $h$ (or a number thereof) \cite{frjg,Chambolle04,Boyd,BoydADMM}. In particular, considering the desired image estimate $x$ to be a (column) vector in $\mathbb{R}^n$, both $f$ and $h$ are usually defined as non-negative functionals on $\mathbb{R}^n$, in which case the standard form of an optimization-based imaging task is given by
\begin{equation}\label{tasks1}
\min_x \,  f(x)+\lambda h(x).
\end{equation}
Here $\lambda > 0$ is a {\em regularization constant} that balances the effects of empirical and prior information on the optimal solution. Specifically, the first (fidelity) term in \eqref{tasks1} forces the solution to ``agree" with the observed data $y$, as it would be the case if one set, e.g., $f(x) = \frac{1}{2} \sum_{i=1}^n |x_i - y_i|^2 = \frac{1}{2} \| x - y \|_2^2$. On the other hand, the prior information is represented by the second (regularization) term in \eqref{tasks1}, which is frequently required to prevent overfitting and/or to render the optimal solution unique and stably computable. For instance, when the optimal solution is expected to be sparse, it is common to set $h(x)=\sum_{i=1}^n |x_i| = \|x\|_1$\cite{frjg,beckteboulle,turlach}.

The convexity and differentiability of the squared Euclidean distance, along with the unparalleled convenience of its numerical handling, are behind the main reasons for its prevalence as a measure of image proximity. The same applies to Mean Squared Error (MSE) and Peak to Signal Noise Ratio (PSNR)---closely related metrics, which have been extensively used to quantify visual quality of images and videos. Yet, neither of the above quantitative measures can be considered a good model for the Human Visual System (HVS), as opposed to the Structural Similarity Index (SSIM), originally proposed by Wang {\em et al.} \cite{wang2002, wang2004}. The SSIM index is based upon the assumption that the HVS has evolved to perceive visual distortions as changes in structural information. On the basis of subjective quality assessments involving large databases, SSIM has been generally accepted to be one of the best measures of visual quality and, by extension, of  perceptual proximity between images.

Considering the exceptional characteristics of the SSIM mentioned above, the idea of replacing the standard norm-based fidelity measures with SSIM seems rather straightforward. However, the optimization problems thus obtained demand much more careful treatment, especially when the question of existence and uniqueness of their solutions is concerned. The main difficulty here stems from the fact that SSIM is not a convex function.

Notwithstanding the above obstacles, optimization problems employing the SSIM as a data fidelity term $f$ have already been addressed. For instance, in \cite{BrunetVrscayWang10}, the authors address the problem of finding the best approximation of data images in the domain of an orthogonal transform, with the optimization performed with respect to the SSIM measure. In this work, instead of maximizing SSIM, the authors minimize
\begin{equation}\label{txy}
T(x,y) = 1-\text{SSIM}(x, y), \quad \mbox{with } x = \Phi c,
\end{equation}
in which case $\Phi$ is a matrix representing the {\em synthesis} phase of an orthogonal transform (e.g., either Fourier, cosine, or wavelet transform), $c$ is a vector of approximation coefficients, and $y$ is a data image to be approximated. It is important to note that, as opposed to the SSIM, $T(x,y)$ may be considered a measure of structural {\em dissimilarity} between $x$ and $y$.

Based on the above results, Rehman {\em et al.} \cite{Rehman} address the problem of sparse dictionary learning through modifying the original k-SVD method of Elad {\em et al.} in \cite{Elad}. The conceptual novelty of the work in \cite{Rehman} consists in using SSIM instead of $\ell_2$-norm based proximity measures, which have been shown to yield substantial improvements in the performance of the sparse learning as well as its application to a number of image processing tasks, including super-resolution. Another interesting use of the SSIM for denoising images was proposed in \cite{Channappayya}. Here, the authors introduce the statistical SSIM index (statSSIM), an extension of the SSIM for wide-sense stationary random processes. Subsequently, using the proposed statSSIM, the authors have been able to reformulate a number of classical statistical estimation problems, such as adaptive filter design. Moreover, it was shown in \cite{Channappayya} that statSSIM is a {\em quasi-concave} function, which opens up a number of possibilities for its efficient numerical optimization, e.g., using the bisection method \cite{Boyd,Channappayya}. Interestingly enough, what seems to have been omitted in \cite{Channappayya} is to recognize that the SSIM is a quasi-concave function \cite{BrunetVrscayWang12} that can be optimized by means of a number of a number of efficient numerical methods as well. For more examples of exploiting SSIM in imaging sciences, the reader is referred to \cite{Wang12,Rehman13,Brunet2017}, which demonstrate the application of this measure to rate distortion optimization, video coding, and image classification.

Finally, we note that maximization of $\text{SSIM}(x,y)$ is equivalent to minimization of $T(x,y)$ in \eqref{txy}, with $T(x, y) = 0$ if and only if $x = y$. Consequently, many image processing problems can be formulated in the form of a quasi-convex program as given by
\begin{eqnarray}\label{eqregssim}
\min_x &&T\left( \Phi(x), y \right) \\
\text{subject to } && h_i(x) \le 0, \quad i = 1, \ldots, p, \notag
\end{eqnarray}
where $x \in\mathbb{R}^n$ is an optimization variable, $\Phi:\mathbb{R}^n\to\mathbb{R}^m$ is a linear transform, $y \in\mathbb{R}^m$ is a vector of observed data, and $h_i:\mathbb{R}^n\to \mathbb{R}$ are convex inequality constraint functions. With little loss of generality, in this work, we assume $p=1$, in which case the problem in \eqref{eqregssim} can be rewritten in its equivalent {\em unconstrained} (Lagrangian) form as
\begin{equation}\label{regssim}
\min_x \, T(\Phi(x),y)+\lambda h(x),
\end{equation}
which fits the format of \eqref{tasks1}. In what follows, we refer to the above problem as an unconstrained SSIM-based optimization problem, as opposed to its constrained version in \eqref{eqregssim}. Note that both problems could also be viewed as special instances under the umbrella of {\em SSIM-based optimization}.

As opposed to developing specific methods for particular applications, in this paper, we introduce a set of algorithms to solve the general problems (\ref{regssim}) and (\ref{eqregssim}). Subsequently, we demonstrate the performance of these algorithms with several applications such as (TV) denoising and sparse approximation, as well as providing comparisons against more traditional approaches.

\section{Structural Similarity Index (SSIM)}

\subsection{Definition}
\label{def}

The SSIM index provides a measure of visual similarity between two images, which could be, for instance, an original scheme and its distorted version. Since it is assumed that the distortionless image is always available, the SSIM is considered a \textit{full-reference} measure of image quality assessment (IQA) \cite{wang2004}. Its definition is based on two assumptions: (i) images are highly structured (that is, pixel values tend to be correlated, especially if they are spatially close), and (ii) HVS is particularly sensitive to structural information. For these reasons, SSIM assesses similarity by quantifying changes in perceived structural information, which can be expressed in terms of the luminance, contrast, and structure of the compared images. In particular, given two images, $x$ and $y$, let $\mu_x$ and $\mu_y$ denote their mean values. Also, let $\sigma_x$ and $\sigma_y$ be the standard deviations of the images, whose cross-correlation coefficient is equal to $\sigma_{x y}$. Then, discrepancies between the luminance of $x$ and $y$ can be quantified using
\begin{equation}
l(x,y)=\frac{2\mu_x\mu_y+C_1}{\mu_x^2+\mu_y^2+C_1},
\end{equation}
where $C_1 > 0$ is added for stability purposes. Note that $l(x,y)$ is sensitive to relative (rather than absolute) changes of luminance, which makes it consistent with Weber's law---a model for light adaptation in the HVS \cite{wang2004}.
\noindent
Further, the contrast of $x$ and $y$ can be compared using
\begin{equation}
c(x,y)=\frac{2\sigma_x\sigma_y+C_2}{\sigma_x^2+\sigma_y^2+C_2},
\end{equation}
where, as before, $C_2 > 0$ is added to prevent division by zero. It is interesting to note that, when there is a change in contrast, $c(x,y)$ is more sensitive if the base contrast is low than when this is high---a characteristic behaviour of the HVS.

\noindent
Finally, the structural component of SSIM is defined as
\begin{equation}
s(x,y)=\frac{\sigma_{xy}+C_3}{\sigma_x\sigma_y+C_3},
\end{equation}
with $C_3 >0$, which makes it very similar to the normalized cross-correlation.

Once the three components $l$, $c$, and $s$ are computed, SSIM is defined according to
\begin{equation}\label{ssim}
\text{SSIM}(x,y)=l(x,y)^\alpha c(x,y)^\beta s(x,y)^\gamma,
\end{equation}
where $\alpha>0$, $\beta>0$ and $\gamma>0$ control the relative influence of their respective components. In \cite{wang2004}, the authors simplify (\ref{ssim}) by setting $\alpha=\beta=\gamma=1$ and $C_3=C_2/2$. This leads to the following well known formula
\begin{equation}
\text{SSIM}(x,y)=\left(\frac{2\mu_x\mu_y+C_1}{\mu_x^2+\mu_y^2+C_1}\right)\left(\frac{2\sigma_{xy}+C_2}{\sigma_x^2+\sigma_y^2+C_2}\right).
\label{stassim}
\end{equation}
This definition of the SSIM will be employed for the remainder of the paper.

We close this section by pointing out that the statistics of natural images vary greatly across their spatial domain, in which case it makes sense to replace the {\em global} means $\mu_x$, $\mu_y$, variances $\sigma_x$, $\sigma_y$, and cross-correlation $\sigma_{x y}$ by their {\em localized} versions computed over $N$ (either distinct or overlapping) local neighbourhoods. In this case, using the localized statistics would result in $N$ values of the SSIM, $\{ {\rm SSIM}_{i}(x_i,y_i)\}_{i=1}^N$, which can be averaged to yield a {\em mean} SSIM index (MSSIM) \cite{wang2004}, which is a more frequently used metric in practice.

\subsection{A normalized metric yielded by the SSIM}

Let $x$ and $y$ be (column) vectors in $\mathbb{R}^n$. In the special case when $x$ and $y$ have equal means, i.e., $\mu_x = \mu_y$, the luminance component $l(x,y)$ in \eqref{stassim} is equal to one, therefore \eqref{stassim} is reduced to
\begin{equation}
\text{SSIM}(x,y)=\frac{2\sigma_{xy}+C_2}{\sigma_x^2+\sigma_y^2+C_2}.
\label{varssim}
\end{equation}
This less cumbersome version of the SSIM can be simplified even further if both $x$ and $y$ have zero mean, i.e., $\mu_x = \mu_y = 0$. In this case, it is rather straightforward to show that
\begin{equation}\label{ssimple}
\text{SSIM}(x,y)=\frac{2x^Ty+C}{\|x\|_2^2+\|y\|_2^2+C}, \quad C = C_2 n,
\end{equation}
with the associated dissimilarity index $T(x,y)$ thus becoming
\begin{equation}\label{dssimple}
T(x,y) = 1-\text{SSIM}(x,y)=\frac{\|x-y\|^2_2}{\|x\|^2_2+\|y\|_2^2+C}.
\end{equation}
Note that, if $C = 0$, $0 \leq T(x,y) \leq 2$, while $T(x,y)=0$ if and only if $x=y$.

$T(x,y)$ given by \eqref{dssimple} is an example of a (squared) normalized metric, which has been discussed in \cite{Brunet, BrunetVrscayWang12}. Thus, while $\text{SSIM}(x,y)$ tells us about how correlated or similar $x$ and $y$ are, $T(x,y)$ gives us a sense of how far $x$ is from a given observation $y$ (in the normalized sense). Moreover, since in the majority of optimization problems the fidelity term quantifies the {\em distance} between a model/estimate and its observation, it seems more suitable for us to proceed with minimization of $T(x,y)$.

We finally note that, throughout the rest of the paper, unless otherwise stated, we will be working with zero-mean images, thus using the definitions of SSIM as given in \eqref{ssimple} and \eqref{dssimple}. This simplification, however, is by no means restrictive. Indeed, many algorithms of modern image processing are implemented under Neumann boundary conditions, thus preserving the mean brightness (or, equivalently, mean value) of input images. Hence, it is rarely a problem to subtract the mean value before the processing starts, while adding it back to the final result.

\subsection{Quasiconvexity and quasiconcavity}
An interesting property of the dissimilarity measure $T$ is that it is quasiconvex over a half-space of $\mathbb{R}^n$. This can be easily proved by using the fact that a function $f:\mathbb{R}^n\to\mathbb{R}$ is quasiconvex if its domain and all its sub-level sets
\begin{equation}
S_\alpha=\{x\in\textbf{dom}~f~|~f(x)\leq\alpha\}, \quad \alpha\in\mathbb{R},
\label{level}
\end{equation}
are convex \cite{Boyd}.
\begin{theorem}
\cite{Otero}
Let $y\in\mathbb{R}^n$ be fixed. Then, $T(x,y)$ is quasiconvex if $x^Ty\geq-\frac{C}{2}$, where $C$ is the stability constant of the dissimilarity measure $T$.
\end{theorem}

Note that, from the relation $T(x,y)-1=-\text{SSIM}(x,y)$, it immediately follows that  ${\rm SSIM}(x,y)$ is a quasiconcave function over the halfspace defined by $x^T y\geq-\frac{C}{2}$. Moreover, using the definition of quasiconcavity, the same line of arguments as in \cite{Otero} can be used to show that  $T(x,y)$ (resp. ${\rm SSIM}(x,y)$) is a quasiconcave (resp. quasiconvex) function, when restricted to the half-space defined by $x^Ty\leq-\frac{C}{2}$.

\section{Constrained SSIM-based Optimization}
A standard quasiconvex optimization problem is defined as follows:
\begin{eqnarray}
\min_x~~&&f(x)\nonumber\\
\label{quasi}
\text{subject to}~~&& h_i(x)\leq 0,~~i=1,\dots,m\\
&&Ax=b\nonumber,
\end{eqnarray}
where $f:\mathbb{R}^n\to\mathbb{R}$ is a quasiconvex cost function, $Ax=b$ (with $A \in \mathbb{R}^{p\times n}$ and $b \in \mathbb{R}^p$) represents a set of $p$ affine equality constraints, and $h_i$ are a set of $m$ convex inequality constraint functions. A standard approach to solving \eqref{quasi} is by means of the {\em bisection method}, which is designed to converge to the desired solution up to a certain accuracy \cite{Boyd}. This method reduces direct minimization of $f$ to a sequence of convex feasibility problems. Using the above approach it is therefore possible to take advantage of the quasiconvexity of $T(x,y)$ to cast SSIM-based optimization problems as quasiconvex optimization problems, which can be subsequently solved by efficient numerical means.

In particular, in what follows, we consider the following quasiconvex optimization problem
\begin{eqnarray}
\min_x~~&& \,  T(\Phi x, y) \nonumber\\
\label{conssim}
\text{subject to}~~&& h_i(x)\leq 0,~~i=1,\dots,m\\
&&Ax=b\nonumber,
\end{eqnarray}
where $\Phi \in \mathbb{R}^{m\times n}$ represents a linear transform and $y \in \mathbb{R}^m$ is a data image. Although the dissimilarity index $T$ above is generally defined as given by \eqref{dssimple}, in practical settings, it is common to deal with non-trivial data images, suggesting $\| y \|_2^2 > 0$. In this case, the stability constant $C$ in \eqref{dssimple} can be set to zero, thereby leading to a simplified expression for $T$ as given by
\begin{equation}
\label{nose}
T(\Phi x ,y)=\frac{\|\Phi x-y\|^2_2}{\|\Phi x\|^2_2+\|y\|_2^2}.
\end{equation}
Note that, for $\|y\|_2 \neq 0$, the above expression is well defined and differentiable for all $x\in\mathbb{R}^n$.

As it was mentioned earlier in this section, the quasiconvex problem \eqref{conssim} can be efficiently solved by means of the bisection algorithm, a particular version of which we introduce next.

\subsection{Quasiconvex optimization}

Before providing details on the proposed optimization procedure, we note that, in image processing, one usually deals with positive-valued vectors $x$. Moreover, it is also frequently the case that the matrix $\Phi$ obeys the property of producing positive-valued results when multiplied by positive-valued vectors. In other words, provided $x \succeq 0$, the above property guarantees that $\Phi x \succeq 0$, which is very typical for situations when $\Phi$ describes the effects of, e.g., blurring and/or subsampling. In such cases, it is reasonable to assume that $(\Phi x)^T y \ge 0$ (provided, of course, $y \succeq 0$ as well), thereby allowing us to consider $T(\Phi x, y)$ as a quasiconvex function of $x$.

Given the latter, in general, a solution of \eqref{conssim} can be found by solving a sequence of feasibility problems \cite{Boyd}. These problems are formulated using a family of convex inequalities, which are defined by means of a family of functions $\phi_\alpha:\mathbb{R}^n\to\mathbb{R}$ such that
\begin{equation}
f(x)\leq\alpha\iff\phi_\alpha(x)\leq0,
\label{iff}
\end{equation}
with the additional property that $\phi_\beta(x)\leq\phi_\alpha(x)$ for all $x$, whenever $\alpha\leq\beta$. In particular, one can show that all the above properties are satisfied by
\begin{equation}
\phi_\alpha(x)=(1-\alpha)\|\Phi(x)-y\|_2^2-2\alpha \, x^T \Phi^T y.
\label{convexfun}
\end{equation}
Consequently, the feasibility problems then assume the following form
\begin{eqnarray}
\text{find}~~&& x\nonumber\\
\label{feasible}
\text{subject to}~~&& \phi_\alpha(x)\leq 0\\
&&h_i(x)\leq 0,~~i=1,\dots,m\nonumber\\
&&Ax=b\nonumber.
\end{eqnarray}
Let $p^\ast$ be the optimal value of \eqref{conssim}. Then, if \eqref{feasible} is feasible, we have that $p^\ast\leq\alpha$, whereas for the infeasible case, $p^\ast>\alpha$.

Using the fact that $0\leq T(\Phi x,y)\leq2$, and defining $\mathbf{1}$ and $\mathbf{0}$ to be vectors in $\mathbb{R}^n$ whose entries are all equal to one and zero, respectively, we propose the following algorithm for solving (\ref{conssim}).
\noindent
\begin{framed}
\begin{center}
{\sc\textbf{\textsc{Algorithm I:}} Bisection method for constrained SSIM-based optimization}
\begin{tabbing}
\textbf{initialize} $x=\mathbf{0}$, $l=0$, $u=2$, $\epsilon>0$;\\
\textbf{data preprocessing $y=y-\frac{1}{n}\mathbf{1}^T y$;}\\
\textbf{while} $u-l>\epsilon$ \textbf{do}\\
\hspace{.2in}$\alpha:=(l+u)/2$;\\
\hspace{.2in}Solve (\ref{feasible});\\
\hspace{.2in}\textbf{if} (\ref{feasible}) is feasible \textbf{then}\\
\hspace{.4in} $u:=\alpha$;\\
\hspace{.2in}\textbf{elseif} $\alpha=1$ \textbf{then}\\
\hspace{.4in}(\ref{conssim}) can not be solved, \textbf{break};\\
\hspace{.2in}\textbf{else}\\
\hspace{.4in}$l:=\alpha$;\\
\hspace{.2in}\textbf{end if}\\
\textbf{end while}\\
\textbf{return} $x$.
\end{tabbing}
\end{center}
\end{framed}
Notice that this method will find a solution $x^\ast$ such that $l\leq f(x^\ast)\leq l+\epsilon$ in exactly $\lceil\log_2((u-l)/\epsilon)\rceil$ iterations \cite{Boyd}, provided that such solution lies in the quasiconvex region of $T(\Phi x,y)$, in which case the algorithm converges to the optimal $p^\ast$ to a precision controllable by the value of $\epsilon$. Once again, we note that the optimal solution $x^\ast$ can be reasonably assumed to lie within the region of quasi-convexity of $T$, since one is normally interested in reconstructions that are positively correlated to the given observation $y$.

\section{Unconstrained SSIM-based Optimization}

In many practical applications, the feasible set of SSIM-based optimization can be described by a single constraint, in which case the optimization problem in \eqref{conssim} can be rewritten in its equivalent unconstrained (Lagrangian) form, that is,
\begin{equation}
\label{unicorn}
\min_{x}\, T(\Phi x, y)+\lambda h(x),
\end{equation}
where $\lambda > 0$ is a regularization parameter and $h: \mathbb{R}^n \to \mathbb{R}$ is a regularization functional, which is often defined to be convex.  As usual, when $\lambda$ is strictly greater than zero, the regularization term is used to force the optimal solution to reside within a predefined ``target" space, thus letting one to use {\em a priori} information on the properties   of the desired solution. This is particularly important in the case where $\Phi^T \Phi$ is either rank-deficient or poorly conditioned, in which case the regularization can help to render the solution well-defined and stably computable. However, since the first term in \eqref{unicorn} is not convex, the entire cost function is not convex either. Consequently, the existence of a unique global minimizer of \eqref{unicorn} cannot be generally guaranteed. With this restriction in mind, however, it is still possible to devise efficient numerical methods capable of converging to either a locally or globally optimal solution, as it will be shown in the following section of the paper.

\subsection{ADDM-based Approach}

In order to solve problem in \eqref{unicorn} we follow an approach based on the Augmented Lagrangian Method of Multipliers (ADMM). This methodology is convenient since it allows us to solve a wide variety of unconstrained SSIM-based optimization problems by splitting the cost function to be minimized into simpler optimization problems that are easier to solve. Moreover, as will be seen in Section \ref{adssimm}, in several cases these simpler problems will have closed form solutions which may be computed very quickly. This is, of course, advantageous since it improves the performance of the resulting algorithm in terms of execution time.

Before going any further, it is worthwhile to mention that one of this simpler optimization problems will usually have the following form,
\begin{equation}
\label{partquad}
\min_x T(\Phi x, y)+\lambda \|x-z\|_2^2,
\end{equation}
where both $y\in\mathbb{R}^m$ and $z\in\mathbb{R}^n$ are given.  We consider this special case separately in the following section.

\subsubsection{Quadratic regularization}
\label{quadratic}

Since the cost function in \eqref{partquad} is twice-differentiable, root-finding algorithms such as the \textit{generalized Newton's method} \cite{Ortega} can be employed to find a local zero-mean solution $x^\ast$.
To such an end, we must first compute the gradient of \eqref{partquad} which satisfies the equation,
\begin{equation}
[s(x^\ast)\Phi^T\Phi+\lambda r(x^\ast)I]x^\ast-\lambda r(x^\ast)z-\Phi^Ty=0.
\end{equation}
Here, $s(x):=1-T(\Phi x,y)$, $r(x):=\|\Phi x\|_2^2+\|y\|_2^2+C$, and $I\in\mathbb{R}^{n\times n}$ is the identity matrix. If we define the following function $f:X\subset\mathbb{R}^n\to\mathbb{R}^n$,
\begin{equation}
\label{thef}
f(x)=[s(x)\Phi^T\Phi+\lambda r(x)I]x-\lambda r(x)z-\Phi^Ty,
\end{equation}
then $x^\ast$ is a (zero-mean) vector in $\mathbb{R}^n$ such that $f(x^\ast)=0$.

Under certain conditions, the convergence of this method is guaranteed
by the \textit{Newton-Kantorovich theorem}, stated below.

\begin{theorem}
\cite{Ortega}
Let $X$ and $Y$ be Banach spaces and $g:X\subset A\to Y$. Assume $g$ is Fr\'{e}chet differentiable on an open convex set $D\subset X$ and
\begin{equation}
\|g'(x)-g'(z)\|\leq K\|x-z\|_2, ~\text{for all $x,z\in D$},
\end{equation}
Also, for some $x_0\in D$, suppose that $g'(x_0)^{-1}$ is defined on all $Y$ and that
\begin{equation}
h:=L\|g'(x_0)^{-1}\|\|g'(x_0)^{-1}g(x_0)\|\leq\frac{1}{2},
\end{equation}
where $\|g'(x_0)^{-1}\|\leq\beta$ and $\|g'(x_0)^{-1}g(x_0)\|\leq\eta$. Set
\begin{equation}
K_1=\frac{1-\sqrt{1-2h}}{K\beta} ~\text{and }K_2=\frac{1+\sqrt{1-2h}}{K\beta},
\end{equation}
and assume that $S:=\{x:\|x-x_0\|\leq K_1\}\subset D$. Then, the Newton iterates
\begin{equation}
x_{k+1}=x_k-g'(x_k)^{-1}g(x_k), ~k\in\mathbb{N},
\end{equation}
are well defined, lie in $S$ and converge to a solution $x^\ast$ of $g(x)=0$, which is unique in $D\cap\{x:\|x-x_0\|\leq K_2\}$. Moreover, if $h<\frac{1}{2}$, the order of convergence is at least quadratic.
\end{theorem}

In our particular case, the Fr\'{e}chet derivative of $f$ is its Jacobian, which we denote as $J_f$,
\begin{equation}
J_f(x)=\Phi^T\Phi x(\nabla s(x))^T+s(x)\Phi^T\Phi+\lambda (x-z)(\nabla r(x))^T+\lambda r(x)I.
\end{equation}
By the previous theorem, convergence is guaranteed in any open convex subset $\Omega$ of $X$, where $X\subset\mathbb{R}^n$, as long as the initial guess $x_0$ satisfies the following condition,
\begin{equation}
L\|J_f(x_0)^{-1}\|\|J_f(x_0)^{-1}f(x_0)\|\leq\frac{1}{2} \, .
\label{inicon}
\end{equation}
Here, $J_f(\cdot)^{-1}$ denotes the inverse of $J_f(\cdot)$, and $L>0$ is a constant that is less or equal than the Lipschitz constant of $J_f(\cdot)$. Indeed, it can be proved that for any open convex subset $\Omega\subset X$, $J_f(\cdot)$ is Lipschitz continuous.

\begin{theorem}
Let $f:X\subset\mathbb{R}^n\to\mathbb{R}$ be defined as in Eq. (\ref{thef}). Then, its Jacobian is Lipschitz continuous on any open convex set $\Omega\subset X$; that is, there exists a constant $L>0$ such that for any $x,w\in\Omega$,
\begin{equation}
\|J_f(x)-J_f(w)\|_F\leq L\|x-w\|_2 \, .
\end{equation}
Here, $\|\cdot\|_F$ denotes the Frobenius norm\footnote{The Frobenius norm of an $m\times n$ matrix $A$ is defined as $\|A\|_F=\sqrt{\sum_{i=1}^m\sum_{j=1}^n|a_{ij}|^2}$.} and
\begin{equation}
L = C_1\|\Phi^T\Phi\|_F+\lambda (C_2\|z\|_2+C_3),~C_1,C_2,C_3>0.
\end{equation}
\begin{proof}
See Appendix.
\end{proof}
\end{theorem}

From this discussion, we propose the following algorithm for solving the problem in \eqref{partquad}.
\begin{framed}
\begin{center}
{\sc\textbf{\textbf{Algorithm II:}} Generalized Newton's method for unconstrained SSIM-based optimization with quadratic regularization}
\begin{tabbing}
\textbf{initialize} Choose $x=x_0$ according to (\ref{inicon});\\
\textbf{data preprocessing} $\bar{y}=\frac{1}{n}\mathbf{1}^Ty$, $y=y-\bar{y}\mathbf{1}$;\\
\textbf{repeat}
\hspace{.2in}$x=x-J_f(x)^{-1}f(x)$;\\
\textbf{until} stopping criterion is satisfied\\
\textbf{return} $x$.
\end{tabbing}
\end{center}
\end{framed}
It is worthwhile to point out that in some cases the calculation of the inverse of $J_f(x)$ may be computationally expensive which, in turn, also makes the $x$-update costly. This difficulty can be addressed by updating the variable $x$
in the following manner,
\begin{equation}
J_f(x^k)x^{k+1}=J_f(x^k)x^k-f(x^k).
\end{equation}
Note that this $x$-update involves the solution of a linear system of equations which can be conveniently done by numerical schemes such as the \textit{Gauss-Seidel} method \cite{BoydADMM}.

\subsubsection{ADMM-based Algorithm}
\label{adssimm}

The problem in \eqref{unicorn} can be solved efficiently by taking advantage of the fact that the objective function is separable. Let us introduce a new variable, $z\in\mathbb{R}^n$. Also, let the independent variable of the regularizing term $h$ be equal to $z$. We then define the following optimization problem,
\begin{eqnarray}
\label{uniadmm}
\min_{x,z}&&T(\Phi x, y)+\lambda h(z),\\ \nonumber
\text{subject to}&& x-z=0.
\end{eqnarray}
Clearly, \eqref{uniadmm} is equivalent to problem \eqref{unicorn}, thus a solution of \eqref{uniadmm} is also a minimizer of the original optimization problem \eqref{unicorn}. Furthermore, notice that \eqref{uniadmm} is presented in ADMM form \cite{BoydADMM}, implying that an analogue of the standard ADMM algorithm can be employed to solve it.

As is customary in the ADMM methodology, let us first form the corresponding augmented Lagrangian of \eqref{uniadmm},
\begin{equation}
\label{Lagran}
L_\rho(x,z,v)=T(\Phi x, y)+\lambda h(z)+v^T(x-z)+\frac{\rho}{2}\|x-z\|_2^2.
\end{equation}
This expression can be rewritten in its more convenient scaled form by combining the linear and quadratic terms and scaling the dual variable $v$, yielding
\begin{equation}
\label{Lagrang}
L_\rho(x,z,u)=T(\Phi x, y)+\lambda h(z)+\frac{\rho}{2}\|x-z+u\|_2^2,
\end{equation}
where $u=v/\rho$. We employ this version of Eq. \eqref{Lagran} since the formulas associated with it are usually simpler than their unscaled counterparts \cite{BoydADMM}. As usual, the iterations of the proposed algorithm for solving \eqref{uniadmm} will be the minimization of Eq. \eqref{Lagrang} with respect to variables $x$ and $z$ in an alternate fashion, and the update of the dual variable $u$, which accounts for the maximization of the dual function $g(u)$:
\begin{equation}
g(u):=\inf_{x,y}L_\rho(x,z,u).
\end{equation}
As such, we define the following iteration schemes for minimizing the cost function of the equivalent counterpart of problem \eqref{unicorn}:
\begin{eqnarray}
x^{k+1}&:=&\argmin_x\left(T(\Phi x,y)+\frac{\rho}{2}\|x-z^{k}+u^{k}\|_2^2\right),\\
z^{k+1}&:=&\argmin_z\left(h(z)+\frac{\rho}{2\lambda}\|x^{k+1}-z+u^{k}\|_2^2\right),\\
u^{k+1}&:=&u^k+x^{k+1}-z^{k+1}.
\end{eqnarray}
Observe that the $x$-update can be computed using the method described in the previous section. Furthermore, when $h$ is convex, the $z$-update is equal to the \textit{proximal operator} of $(\lambda/\rho)h$ \cite{Parikh}. Recall that for a convex function $f:\mathbb{R}^n\to\mathbb{R}$ its proximal operator $\mathbf{prox}_f:\mathbb{R}^n\to\mathbb{R}^n$ is defined as
\begin{equation}
\mathbf{prox}_f(v):=\argmin_x\left(f(x)+\frac{1}{2}\|x-v\|_2^2\right) \, .
\end{equation}
It follows that
\begin{equation}
z^{k+1}:=\mathbf{prox}_{\frac{\lambda}{\rho}h}(x^{k+1}+u^k).
\end{equation}

Taking this into account, we introduce the following algorithm for solving the problem in \eqref{unicorn}.
\begin{framed}
\begin{center}
{\sc\textbf{\textsc{Algorithm III:}} ADMM-based method for unconstrained SSIM-based optimization}
\begin{tabbing}
\textbf{initialize} $x=z=x_0$, $u=\mathbf{0}$;\\
\textbf{data preprocessing} $y=y-\frac{1}{n}\mathbf{1}^T y$;\\
\textbf{repeat}\\
\hspace{.2in}$x:=\argmin_x\left(T(\Phi x,y)+\frac{\rho}{2}\|x-z+u\|_2^2\right)$;\\
\hspace{.2in}$z:=\argmin_z\left(h(z)+\frac{\rho}{2\lambda}\|x-z+u\|_2^2\right)$;\\
\hspace{.2in}$u:=u+x-z$;\\
\textbf{until} stopping criterion is satisfied.\\
\textbf{return} $x$.
\end{tabbing}
\end{center}
\end{framed}

\subsection{Alternative ADMM-based Approach}

As shown in the previous section, one of the advantages of the ADMM-based approach is that it allows us to solve problem \eqref{unicorn} by solving a sequence of simpler optimization problems.  The complexity of the proposed method can be further reduced by reformulating \eqref{unicorn} as the following equivalent optimization problem:
\begin{eqnarray}
\min_{x,w,z}&&T(w, y)+\lambda h(z),\\ \nonumber
\text{subject to}&& \Phi x-w=0, \nonumber\\
&& x-z=0.
\label{altadmm}
\end{eqnarray}
The scaled form of the augmented Lagrangian of this new problem is given by
\begin{eqnarray*}
L_{\rho,\mu}(x,w,z,u,v)&=&T(w, y)+\lambda h(z)+\frac{\rho}{2}\|\Phi x-w+u\|_2^2+\frac{\mu}{2}\|x-z+v\|_2^2\nonumber.
\end{eqnarray*}
Therefore, problem \eqref{altadmm} can be solved by defining the following iterations:
\begin{eqnarray}
x^{k+1}&:=&\argmin_x\left(\frac{\rho}{2}\|\Phi x-w^k+u^{k}\|_2^2+\frac{\mu}{2}\|x-z^{k}+v^{k}\|_2^2\right),\\
w^{k+1}&:=&\argmin_w\left(T(w,y)+\frac{\rho}{2}\|w-\Phi x^{k+1}-u^{k}\|_2^2\right),\label{wup}\\
z^{k+1}&:=&\argmin_z\left(h(z)+\frac{\mu}{2\lambda}\|z-x^{k+1}-v^{k}\|_2^2\right),\\
u^{k+1}&:=&u^k+\Phi x^{k+1}-w^{k+1},\\
v^{k+1}&:=&v^k+ x^{k+1}-z^{k+1}.
\end{eqnarray}

\noindent
Notice that the $x$-update has a closed form solution, which is given by
\begin{equation}
x^{k+1}=(\rho\Phi^T\Phi+\mu I)^{-1}(\rho\Phi^T(w^k-u^k)+\mu(z^k-v^k)),
\end{equation}
where $I\in\mathbb{R}^{n\times n}$ is the identity matrix. Alternatively, observe that the variable $x$ can also be updated by solving the following system of linear equations:
\begin{equation}
(\rho\Phi^T\Phi+\mu I)x^{k+1}=(\rho\Phi^T(w^k-u^k)+\mu(z^k-v^k)).
\end{equation}
This option is convenient when the computation of the inverse of the matrix $\rho\Phi^T\Phi+\mu I$ is expensive. Since the matrix $\rho\Phi^T\Phi+\mu I$ is symmetric and positive-definite we can employ a \textit{conjugate gradient method} to update the primal variable $x$ \cite{Hestenes}.

The $w$-update is a simple and special case of the more general method described in Section \ref{quadratic}, thus the methods discussed in that section can be employed to compute \eqref{wup}.

Finally, the $z$-update is simply the proximal operator of the function $\frac{\lambda}{\mu}h:\mathbb{R}^n\to\mathbb{R}$,
\begin{equation}
z^{k+1}:=\mathbf{prox}_{\frac{\lambda}{\mu}h}(x^{k+1}+u^k),
\end{equation}
provided $h$ is a convex function of $x$.

Given the above discussion, we present the following alternative algorithm for solving the problem in \eqref{unicorn}.
\newpage
\begin{framed}
\begin{center}
{\sc\textbf{\textsc{Algorithm IV:}} Alternative ADMM-based method for unconstrained SSIM-based optimization}
\begin{tabbing}
\textbf{initialize} $x=z=x_0$, $w=\Phi x_0$, $u=v=\mathbf{0}$;\\
\textbf{data preprocessing} $y=y-\frac{1}{n}\mathbf{1}^T y$;\\
\textbf{repeat}\\
\hspace{.2in}$x:=(\rho\Phi^T\Phi+\mu I)^{-1}(\rho\Phi^T(w-u)+\mu(z-v))$;\\
\hspace{.2in}$w:=\argmin_w\left(T(w,y)+\frac{\rho}{2}\|w-\Phi x-u\|_2^2\right)$;\\
\hspace{.2in}$z:=\argmin_z\left(h(z)+\frac{\mu}{2\lambda}\|z-x-v\|_2^2\right)$;\\
\hspace{.2in}$u:=u+\Phi x-w$;\\
\hspace{.2in}$v:=v+x-z$;\\
\textbf{until} stopping criterion is satisfied.\\
\textbf{return} $x$.
\end{tabbing}
\end{center}
\end{framed}

\section{Applications}

Naturally, different sets of constraints and regularization terms yield different SSIM-based imaging problems. In this section, we review some applications and the corresponding optimization problems associated with them.

\subsection{SSIM with Tikhonov regularization}

A common method used for ill-posed problems is \textit{Tikhonov regularization} or \textit{ridge regression}. This is basically a constrained version of least squares and it is found in different fields such as statistics and engineering. It is stated as follows
\begin{eqnarray}
\label{tiko}
\min_x&&\,\|\Phi x-y\|_2^2\\
\text{subject to}&&\, \|Ax\|_2^2\leq \lambda\nonumber,
\end{eqnarray}
where $A\in\mathbb{R}^{p\times n}$ is called the \textit{Tikhonov matrix}. A common choice for the matrix $A$ is the identity matrix, however, other choices may be a scaled finite approximation of  a differential operator or a scaled orthogonal projection \cite{HochstenbachReichel,FuhryReichel,MorigiReichelSgallari}.

The SSIM counterpart of problem (\ref{tiko}) is obtained by replacing the Euclidean fidelity term of problem (\ref{tiko}) by the dissimilarity measure $T(\Phi x,y)$:
\begin{eqnarray}
\min_x&&\,\,T(\Phi x,y)\\
\text{subject to}&&\,\|Ax\|_2^2\leq \lambda\nonumber.
\end{eqnarray}
Notice that this problem can be addressed by employing Algorithm I.

Algorithm II can be employed to solve the unconstrained counterpart of the previous problem, given by
\begin{equation}
\min_{x} \, T(\Phi x,y)+\lambda\|Ax\|_2^2 \, .
\end{equation}
In this particular case, the corresponding iteration schemes
are defined as follows,
\begin{eqnarray}
x^{k+1}&:=&\argmin_x\left(T(\Phi x,y)+\frac{\rho}{2}\|x-z^{k}+u^{k}\|_2^2\right),\\
z^{k+1}&:=&\left(\frac{2\lambda}{\rho}A^TA+I\right)^{-1}(x^{k+1}+u^{k}),\\
u^{k+1}&:=&u^k+x^{k+1}-z^{k+1},
\end{eqnarray}
where $I\in\mathbb{R}^{n\times n}$ is the identity matrix.  In some cases, the computation of the inverse of the matrix $\frac{2\lambda}{\rho}A^TA+I$ may be expensive, making it more appropriate to employ alternative methods to compute the $z$-update. For instance, if $p$ is smaller than $n$, it may be convenient to employ the \textit{matrix inversion lemma} \cite{BoydADMM}. Furthermore, notice that the $z$-update is equivalent to solving the following system of linear equations,
\begin{equation}
\left(\frac{2\lambda}{\rho}A^TA+I\right)z^{k+1}:=x^{k+1}+u^{k} \, .
\end{equation}
Since the matrix $\frac{2\lambda}{\rho}A^TA+I$ is positive-definite and symmetric, the variable $z$ can be updated efficiently by using a conjugate gradient method.

\subsection{SSIM-$\mathbf{\ell^1}$-based optimization}
\label{eleuno}

Other interesting applications emerge when the $\ell_1$ norm is used either as a constraint or a regularization term. For example, by defining the constraint $h(x)=\|x\|_1-\lambda$, we obtain the following SSIM-based optimization problem,
\begin{eqnarray}
\label{ssiml1}
\min_x&& \, T(\Phi x,y)    \\
\text{subject to}&& \|x\|_1\leq \lambda\nonumber.
\end{eqnarray}
As expected, the solution of \eqref{ssiml1} can be obtained by means of Algorithm I. Clearly, the unconstrained counterpart of (\ref{ssiml1}) is given by,
\begin{equation}
\min_{x} \, T(\Phi x,y)+\lambda\|x\|_1,
\label{unssiml1}
\end{equation}
which can be solved with Algorithm II. Problem \eqref{unssiml1} can be seen as the SSIM version of the classical regularized version of the least squares method known as LASSO (Least Absolute Shrinkage and Selection Operator) \cite{efron2004,frjg}. The corresponding iteration schemes to solve it are as follows,
\begin{eqnarray}
x^{k+1}&:=&\argmin_x\left(T(\Phi x,y)+\frac{\rho}{2}\|x-z^{k}+u^{k}\|_2^2\right),\\
z^{k+1}&:=&S_{\frac{\lambda}{\rho}}(x^{k+1}+u^k),\\
u^{k+1}&:=&u^k+x^{k+1}-z^{k+1} \, .
\end{eqnarray}
Here, the $z$-update is equal to the proximal operator of the $\ell_1$ norm, also known as the \textit{soft thresholding operator} \cite{efron2004,frjg,BoydADMM}. This is an element-wise operator, which is defined as
\begin{equation}
  S_\tau(t)=
  \begin{cases}
  t-\tau, & t>\tau,\\
  0, & |t|\leq\tau,\\
  t+\tau, & t<\tau\, .
  \end{cases}
\end{equation}

Problems \eqref{ssiml1} and \eqref{unssiml1} are appealing because they combine the concepts of similarity and sparseness, therefore, these are relevant in applications in which sparsity is the assumed underlying model for images. To the best of our knowledge, optimization of the SSIM along with the $\ell^1$ norm has been only reported in \cite{Otero,Otero2014,OteroIPCV} and in this work.

\subsection{SSIM and Total Variation}
\label{latv}

By employing the constraint $h(x)=\|Dx\|_1-\lambda$, where $D\in\mathbb{R}^{n\times n}$ is a difference matrix and $\Phi$ is the identity matrix $I\in\mathbb{R}^{n\times n}$, we can define a SSIM-total-variation-denoising method for one-dimensional discrete signals as follows.  Given a noisy signal $y\in\mathbb{R}^n$, its denoised reconstruction $x$ is the solution of the problem,
\begin{eqnarray}
\label{consTV}
\min_x&& \,\, T (x,y)   \\
\text{subject to}\,&& \|Dx\|_1\leq \lambda\nonumber.
\end{eqnarray}
Here, we consider $\|Dx\|_1$ as a measure of the total variation (TV) of the vector $x$. Notice that instead of minimizing the TV, we employ it as a constraint. Moreover, methods for solving the $\ell_2$ version of  \eqref{consTV} can be found in \cite{CombettesP04,FadiliP11}. As with the classical TV optimization problems, solutions of \eqref{consTV} have bounded variation as well.

As for images, these can be denoised in the following manner. Let $Y\in\mathbb{R}^{m\times n}$ be a noisy image. Also, let $V:\mathbb{R}^{m\times n}\to\mathbb{R}^{mn\times1}$ be a linear transformation that converts matrices into column vectors, that is,
\begin{equation}
V(A)=\vect(A)=[a_{11}, a_{21},\dots,a_{(m-1)n},a_{mn}]^T.
\end{equation}
A reconstruction $X\in\mathbb{R}^{m\times n}$ of the noiseless image can be obtained by means of the following SSIM-based optimization problem,
\begin{equation}
\min_{X} \, T(V(X),V(Y))+\lambda\|X\|_{TV},
\label{unssimtv}
\end{equation}
where the regularizing term is a discretization of the isotropic TV seminorm for real-valued images. In particular, the TV of $X$ is the discrete counterpart of the standard total variation of a continuous image $g$ that belongs to the function space $L^1(\Omega)$, where $\Omega$ is an open subset of $\mathbb{R}^2$:
\begin{equation}
TV(g)=\int_\Omega \|Dg(x)\|_2dx=\sup_{\xi\in \Xi}\left\{\int_\Omega g(x)\nabla\cdot\xi~dx\right\}.
\label{dualtv}
\end{equation}
Here, $\Xi=\{\xi:\xi\in C_c^1(\Omega,\mathbb{R}^2), \|\xi_k(x)\|_2\leq 1\ \forall x\in\Omega\}$, and $\nabla\cdot$ is the divergence operator \cite{Chambolle04}. As anticipated, we employ the following iterations for solving \eqref{unssimtv}:
\begin{eqnarray}
X^{k+1}&:=&\argmin_X\left(T(V(X),V(Y))+\frac{\rho}{2}\|X-Z^{k}+U^{k}\|_F^2\right),\\
Z^{k+1}&:=&\argmin_Z\left(\|Z\|_{TV}+\frac{\rho}{2\lambda}\|Z-X^{k+1}-U^{k}\|_F^2\right),\\
U^{k+1}&:=&U^k+X^{k+1}-Z^{k+1},
\end{eqnarray}
where $\|\cdot\|_F$ is the Frobenius norm. Notice that the $Z$-update may be computed efficiently by using the algorithm introduced by Chambolle in \cite{Chambolle04}.

As mentioned before, it is more convenient to employ an average of local SSIMs as a fidelity term. Let $\{Y_i\}_{i=1}^N$ be a partition of the given image $Y$ such that $\cup_{i=1}^NY_i=Y$. Further, let $\{X_i,Z_i\}_{i=1}^N$ also be partitions of the variables $X$ and $Z$ such that $\cup_{i=1}^NX_i=X$ and $\cup_{i=1}^NZ_i=Z$. Also, let $MT:\mathbb{R}^{m\times n}\times\mathbb{R}^{m\times n}\to\mathbb{R}$ be given by
\begin{equation}
MT(X,Y)=\frac{1}{N}\sum_{i=1}^N T(V(X_i),V(Y_i)).
\label{MT}
\end{equation}
Then, the optimization problem that is to be solved is
\begin{equation}
\min_{X} \, MT(X,Y)+\lambda\|X\|_{TV}.
\label{unmssimtv}
\end{equation}
If $\{Y_i,X_i\,Z_i\}_{i=1}^N$ are partitions of non-overlapping blocks, the problem in \eqref{unmssimtv} can be solved by carrying out the following iterations,
\begin{eqnarray}
X_i^{k+1}&:=&\argmin_{X_i}\left(T(V(X_i),V(Y_i))+\frac{N\rho}{2}\|X_i-Z_i^{k}+U_i^{k}\|_F^2\right),\\
Z^{k+1}&:=&\argmin_Z\left(\|Z\|_{TV}+\frac{\rho}{2\lambda}\|Z-X^{k+1}-U^{k}\|_F^2\right),\\
U^{k+1}&:=&U^k+X^{k+1}-Z^{k+1},
\end{eqnarray}
where $U_i$ is an element of the partition of the dual variable $U$. As expected, $\cup_{i=1}^NU_i=U$, and $U_i\cap U_j=\varnothing$ for all $i\neq j$. The extension of this algorithm when a weighted average of local SSIMs is used as a measure of similarity between images is straightforward.

Other imaging tasks can be performed by solving the following variant of the SSIM-based optimization problem in \eqref{unmssimtv},
\begin{equation}
\min_{X} \, MT(A(X),Y)+\lambda\|X\|_{TV},
\label{unmopssimtv}
\end{equation}
where $A(\cdot)$ is a linear operator (e.g., blurring kernel, subsampling procedure, etc.). A minimizer of \eqref{unmopssimtv} may be found by means of Algorithm IV. In this case, the corresponding ADMM iterations are the following:
\begin{eqnarray}
X^{k+1}&:=&\argmin_X\left(\frac{\rho}{2}\|A(X)-W^k+U^{k}\|_2^2+\frac{\mu}{2}\|X-Z^{k}+V^{k}\|_2^2\right),\\
W^{k+1}&:=&\argmin_W\left(MT(W,Y)+\frac{\rho}{2}\|W-A(X)^{k+1}-U^{k}\|_2^2\right),\\
Z^{k+1}&:=&\argmin_Z\left(\|Z\|_{TV}+\frac{\mu}{2\lambda}\|Z-X^{k+1}-V^{k}\|_2^2\right),\\
U^{k+1}&:=&U^k+A(X)^{k+1}-W^{k+1},\\
V^{k+1}&:=&V^k+ X^{k+1}-Z^{k+1}.
\end{eqnarray}

Furthermore, under certain circumstances, Algorithm III can be employed to solve problem \eqref{unmopssimtv}. Let $\{A(X)_i\}_{i=1}^N$ be a partition of $A(X)$ such that $\cup_{i=1}^NA(X)_i=A(X)$, and $A(X)_i\cap A(X)_j=\varnothing$ for all $i\neq j$. If there exist operators $D_i$, $1\le i\leq N$, such that
\begin{equation}
D_i(V(X_i))=V(A(X)_i),
\end{equation}
for all $i\in\{1,\dots,N\}$, then a minimizer of \eqref{unmopssimtv} can be found by means of Algorithm III. The corresponding iterations are as follows:
\begin{eqnarray}
X_i^{k+1}&:=&\argmin_{X_i}\left(T(D_i(V(X_i)),V(Y_i))+\frac{N\rho}{2}\|X_i-Z_i^{k}+U_i^{k}\|_F^2\right),\\
Z^{k+1}&:=&\argmin_Z\left(\|Z\|_{TV}+\frac{\rho}{2\lambda}\|Z-A(X)^{k+1}-U^{k}\|_F^2\right),\\
U^{k+1}&:=&U^k+X^{k+1}-Z^{k+1},
\end{eqnarray}
where $X_i$, $U_i$ and $Z_i$ are elements of the partitions $\{X_i,U_i,Z_i\}$ defined above.

We close this section by mentioning that to the best of our knowledge, the contributions reported in \cite{YuShao,OteroIPCV,Otero} along with the applications presented above are the only approaches in the literature that combine TV and the SSIM.

\section{Experiments}

In this section, we provide some numerical and visual results to evaluate the performance of some of the methods introduced above. We have focussed our attention on two types of unconstrained SSIM-based optimization problems that have been barely explored in the literature: SSIM with $\ell_1$ regularization and SSIM with the TV seminorm as a regularizing term. The results are presented as follows:  In Section \ref{unconop}, we assess the efficacy of the ADMM-SSIM-based approach when the $\ell_1$ regularization is employed, whereas in Section \ref{ssimtvexp}, performance of the ADMM-SSIM-based methodology that employs the TV seminorm as a regularizing term is evaluated.

In all experiments, we employed non-overlapping pixel blocks. Performance of the $\ell_2$- and SSIM-based approaches is compared by computing the MSSIM of the original images and their corresponding reconstructions. Here, the MSSIM is simply the average of the SSIM values of all non-overlapping blocks, which are computed using Eq. \eqref{nose}.

\newpage

\subsection{SSIM with $\ell_1$ Regularization}
\label{unconop}

Here, we consider the unconstrained SSIM-based optimization problem,
\begin{equation}
\label{missl1}
\min_x \,  T(Dx,y) + \lambda \|x\|_1,
\end{equation}
where $D$ is a $n\times n$ Discrete Cosine Transform (DCT) matrix, and $y\in\mathbb{R}^n$ is the given observation. This problem can be solved with either Algorithm III or Algorithm IV. For these experiments we have employed Algorithm III.

We shall compare the solutions obtained by the proposed method with those obtained by the $\ell_2$ version of problem \eqref{missl1}, namely,
\begin{equation}
\label{eltwo}
\min_x \, \frac{1}{2} \| Dx - y \|_2^2 + \lambda \| x \|_1,
\end{equation}
which can be solved by means of the soft thresholding (ST) operator \cite{frjg,turlach} if $D$ is an orthogonal matrix. For the sake of a fair comparison between the two approaches, the regularization of each $8\times8$ image-block was adjusted so that the $\ell_0$ norm of the set of recovered DCT coefficients is 18 in all cases.

In these experiments, the well-known test images \textit{Lena} and \textit{Mandrill} were employed.


The SSIM maps clearly show that the proposed method (ADMM-SSIM) outperforms the classical $\ell_2$ approach (ST) with respect to the SSIM. Moreover, by taking a closer look at the reconstructions, we observe that the SSIM-based approach yields images which are brighter and posses a higher contrast than the $\ell_2$ reconstructions (e.g., compare the central regions of the reconstructions in the bottom row). This should not be surprising since the dissimilarity measure $T(Dx,y)$ takes into account the component of the SSIM that measures the contribution of the contrast of the images that are being compared. On the other hand, thanks to the structural component of the SSIM, textures are better captured by the proposed method. 

Regarding numerical results, the MSSIMs that are obtained with the proposed approach are 0.9164 (\textit{Lena}) and 0.8440 (\textit{Mandrill}). As for the classical $\ell_2$ method, the corresponding MSSIMs are 0.9020 (\textit{Lena}) and 0.7935 (\textit{Mandrill}). Even though in the case of the image \textit{Lena} the performance of the SSIM-based approach is somewhat superior, the proposed method significantly outperforms the $\ell_2$ reconstruction of \textit{Mandrill} (see the SSIM maps of the third row). The main reason is that this image is much less regular than \textit{Lena}, so its features are better approximated by SSIM-based techniques.

\begin{figure}[htbp]
\begin{center}
\includegraphics[width=1\textwidth,height=1\textwidth]{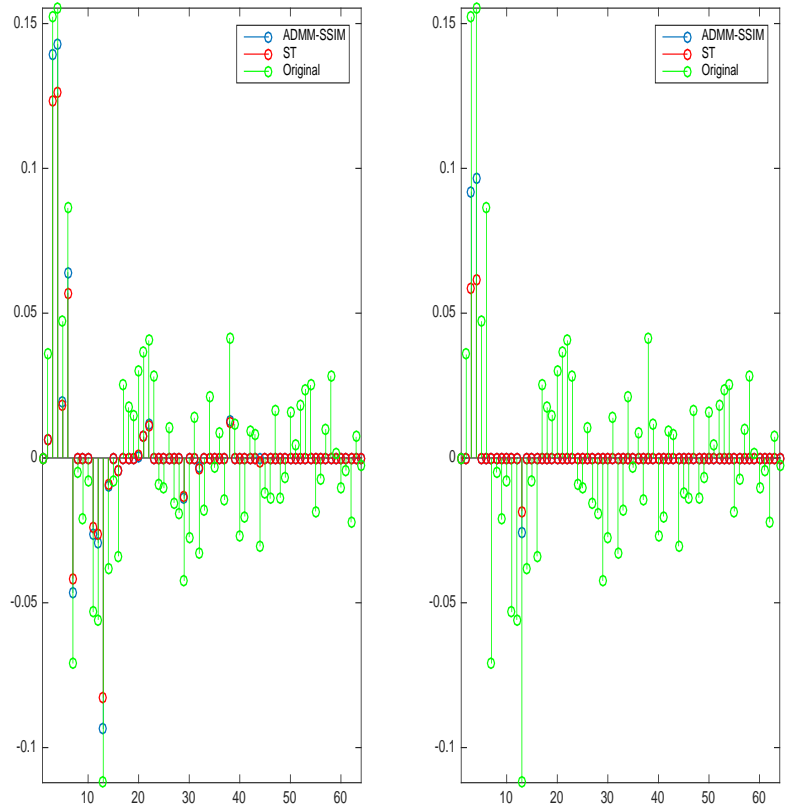}
\caption{A visual comparison between the original and recovered coefficients from a particular block of the \textit{Lena} image can be observed. Regularization is carried out so that the two methods being compared induce the same sparseness in their recoveries. In the two shown examples, the same block was processed but subjected to different amounts of regularization. In particular, the $\ell_0$ norm of the set of DCT coefficients that were recovered by both the proposed method and ST is 18 for the first example (first plot from left to right), and 3 for the second (plot on the right hand side).}
\label{dctcoeff}
\end{center}
\end{figure}

In order to obtain a deeper insight as to how the SSIM-based approach yields brighter reconstructions with higher contrast, it is worthwhile to see the type of solutions that are obtained by both $\ell_2$ and SSIM approaches when the regularization varies. Experimental results show that, in general, the shrinking of the DCT coefficients is not as strong as the classical ST operator. In other words, the SSIM solution is usually a ``scaled'' version of the $\ell_2$ solution, nevertheless, the experimental results show that this scaling is not always the same for all recovered coefficients---an example is presented in Figure \ref{dctcoeff}. In the plots that are shown, the same image-block is processed but subjected to different degrees of regularization. In the plot on the left hand side, the $\ell_0$ norm of both solutions is 18, whereas on the plot on the right hand side, the number of non-zero coefficients for both methods is 3.

\begin{figure}[htbp]
\begin{center}
\includegraphics[width=1\textwidth,height=1\textwidth]{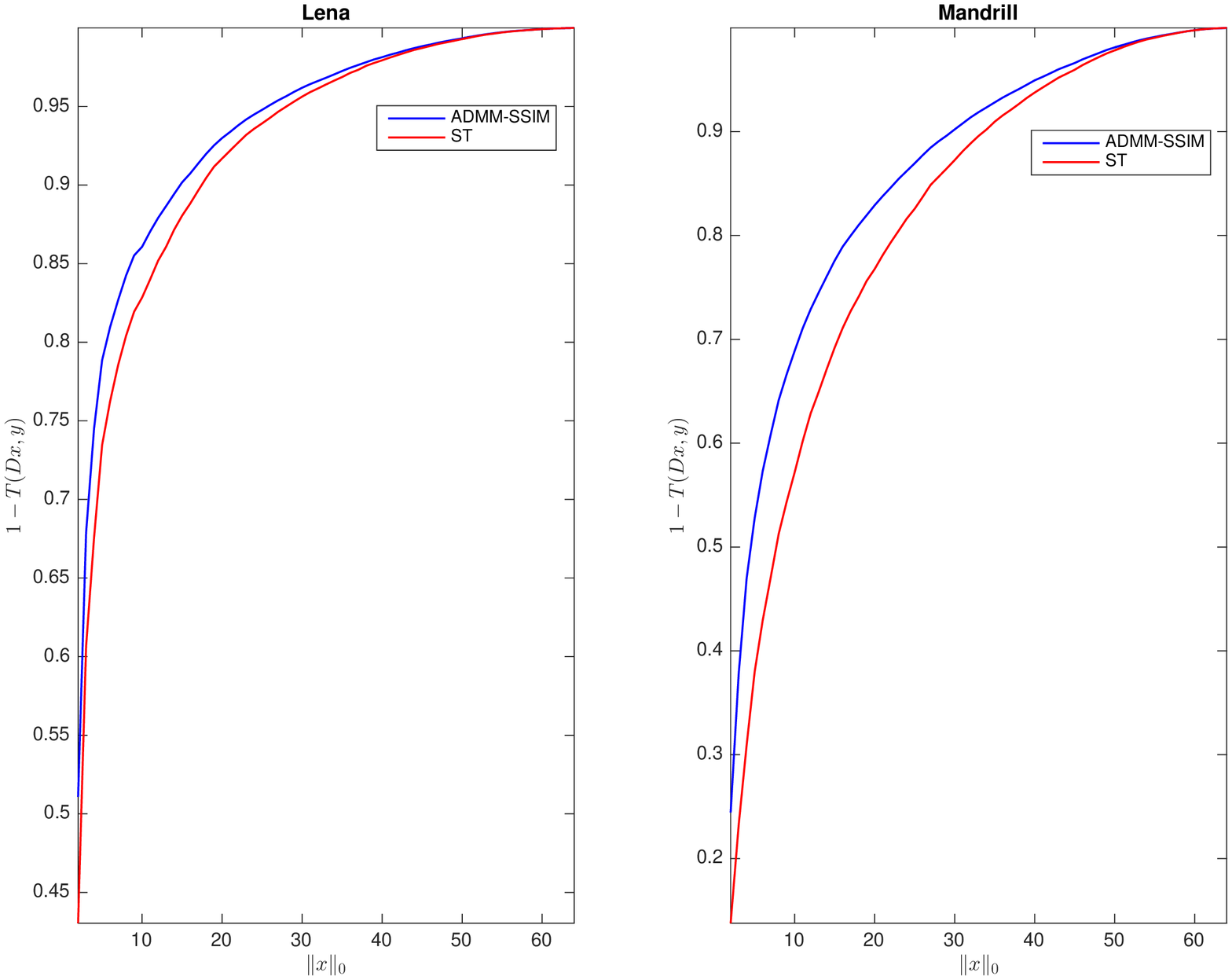}
\caption{In this figure, both plots correspond to the average SSIM versus the $\ell_0$ norm of the recovered coefficients for the test images \textit{Lena} and \textit{Mandrill}. It can be observed that the SSIM-based technique gradually outperforms the classical $\ell_2$ method as regularization increases.}
\label{spassim}
\end{center}
\end{figure}

As a complement to the above discussion, Figure \ref{spassim} shows how the MSSIM changes as a function of the $\ell_0$ norm of the solutions that are obtained by both methods. The plot on the left hand side shows this behaviour of the MSSIM for an image patch of \textit{Lena}, whereas the other plot illustrates what happens in the case of an image patch of \textit{Mandrill}. It can be seen that when the regularization is not strong, the performance of both approaches is very similar.  However, as the regularization is increased, the difference of the MSSIM values yielded by both methods tends to be greater. As expected, this is more noticeable for \textit{Mandrill}.  These results show that when high compression is required, and when the images are not so regular, it is more convenient to opt for a SSIM-based technique for the sake of visual quality.

\subsection{SSIM and Total Variation}
\label{ssimtvexp}

In this section, we examine several imaging tasks which employ the TV seminorm as a regularization term.  In order to assess the performance of the proposed ADMM-SSIM methods, we shall compare their results with their $\ell_2$ counterparts.

It is important to mention that in order to reduce blockiness in the reconstructions the mean of each non-overlapping pixel block is not subtracted prior to processing. This implies that the fidelity term defined in \eqref{MT} is not equivalent but is based on the dissimilarity measure introduced in Section \ref{def}. Despite this, the experiments presented below suggest that this fidelity measure may be used as a substitute of the SSIM.

\subsubsection{Denoising}

In the following experiments, the denoising of some images corrupted with Additive White Gaussian Noise (AWGN) was performed. Although from a maximum a posteriori (MAP) perspective the ADMM-SSIM approach is not optimal, it is worthwhile to see how denoising is carried out when the SSIM-based metric is employed as a fidelity term.

As one might expect, the noiseless approximation is obtained by solving the following unconstrained SSIM-based optimization problem,
\begin{equation}
\min_{X} \, MT(X,Y)+\lambda\|X\|_{TV},
\label{STV}
\end{equation}
where $MT:\mathbb{R}^{m\times n}\times\mathbb{R}^{m\times n}\to\mathbb{R}$ is the fidelity term that was previously defined in Eq. \eqref{MT}. Problem \eqref{STV} can be solved by using either Algorithm III or Algorithm IV.  Here, we employed Algorithm III since the ADMM-SSIM iterations are quite simple optimization problems.

In order to assess the performance of the proposed ADMM-SSIM method, we compare it with its $\ell_2$ counterpart, namely,
\begin{equation}
\min_{X} \, \|X-Y\|_2^2+\lambda\|X\|_{TV}.
\end{equation}
Naturally, Chambolle's algorithm can be employed for solving this optimization problem \cite{Chambolle04}. In order to compare the effectiveness of the proposed approach and Chambolle's method (TV), regularization was carried out in such a way that the TV seminorms of the reconstructions yielded by both methods are the same.

\begin{figure}[htbp]
\centering
\includegraphics[width=1\textwidth,height=0.63\textwidth]{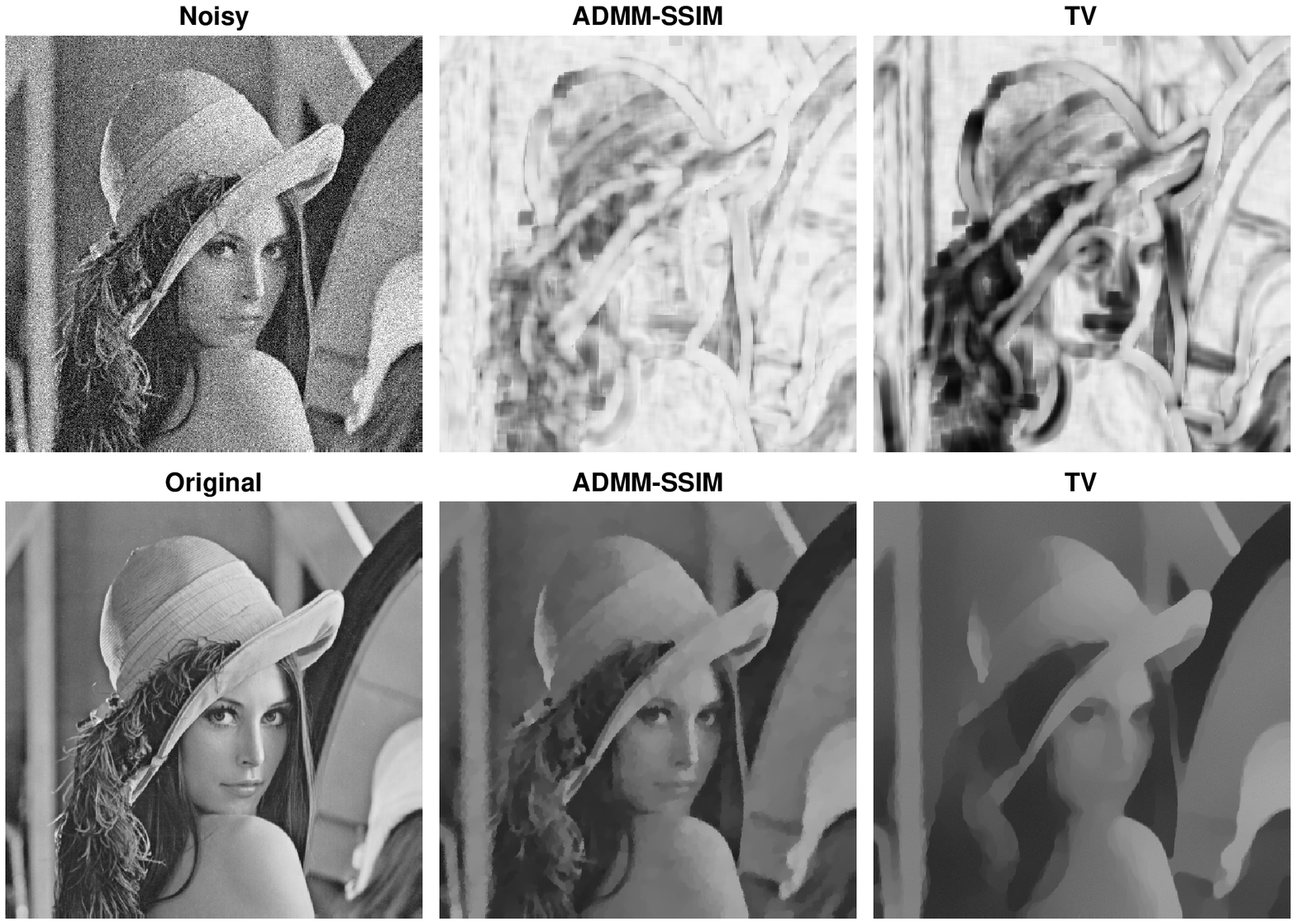}
\includegraphics[width=1\textwidth,height=0.63\textwidth]{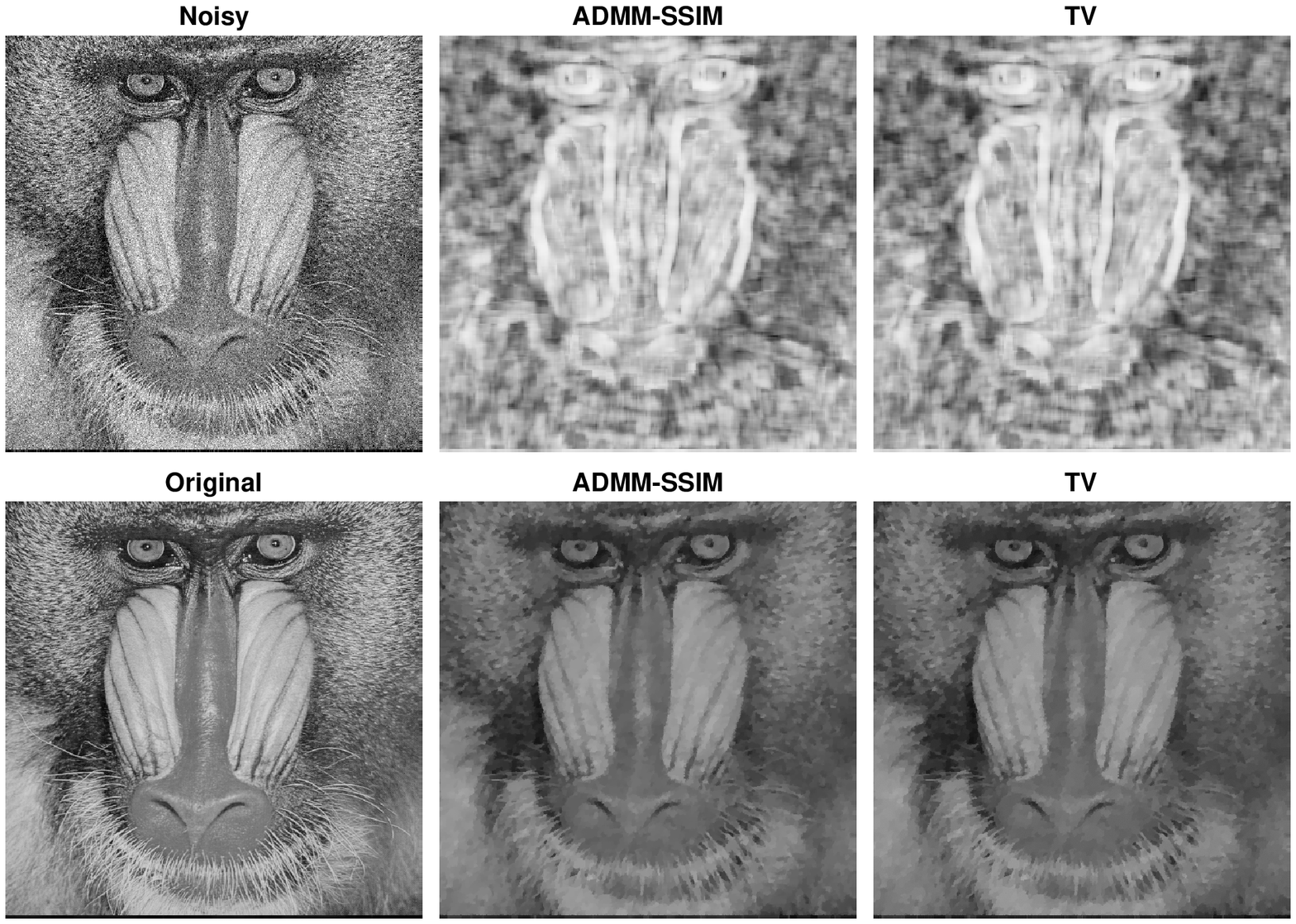}
\caption{Some visual results for the denoising of the test images \textit{Lena} and \textit{Mandrill}. For \textit{Lena}, the TV seminorm of the shown reconstructions is 2500, whereas for \textit{Mandrill} it is 4500. In the first and third rows are shown the noisy images along with the the SSIM maps between the reconstructions and the original images. The original and reconstructed images are shown in the second and fourth rows. The MSSIMs for \textit{Lena} and \textit{Mandrill} that are obtained by using the ADMM-SSIM-based method are 0.6426 and 0.4832, respectively. The corresponding MSSIM values for the
$\ell^2$ approach are 0.4669 (\textit{Lena}) and 0.4708 (\textit{Mandrill}).}
\label{lemanssimtv}
\end{figure}

In Figure \ref{lemanssimtv}, some visual results are shown. Once again, we employed the test images \textit{Lena} and \textit{Mandrill}. The noisy images, as well as the SSIM maps, can be observed in the first and the third rows. Reconstructions and original images are presented in the second and fourth rows. The TV seminorm of the reconstructions is 2500 for \textit{Lena} and 4500 for \textit{Mandrill}. The Peak Signal-to-Noise Ratio (PSNR) prior to denoising was 18.067 dB in all experiments.

In the case of \textit{Lena}, it is evident that the proposed method performs significantly better than its $\ell_2$ counterpart. Notice that some features of the original \textit{Lena} are better reconstructed (e.g., the eyes in \textit{Lena}) whereas in the TV reconstruction these features are considerably blurred. This is mainly due to the fact that the noise does not completely hide some of the more important attributes of the original image. Since the fidelity term enforces the minimizer of problem \eqref{STV} to be visually as similar as possible as the given noisy observation, while denoising is still accomplished with the regularizing term, the reconstruction yielded by the ADMM-SSIM approach is visually more similar to the noiseless image. As for MSSIM values, these are 0.4669 and 0.6426 for the TV and ADMM-SSIM reconstructions, respectively.

Nevertheless, in some circumstances, there is not such a noticeable gap between the ADMM-SSIM and TV approaches.  This is the case for the \textit{Mandrill} image.  This image is much less regular than \textit{Lena}, therefore TV-based denoising techniques will deliver reconstructions devoid of the fine details that the original \textit{Mandrill} has (e.g., the fur).  The ADMM-SSIM method performs better than the $\ell_2$ approach.  However, its effectiveness is only somewhat better due to the low regularity of \textit{Mandrill}. Regarding numerical results, the computed MSSIMs are 0.4832 (ADMM-SSIM) and 0.4708 (TV).

\begin{figure}[htbp]
\begin{center}
\includegraphics[width=1\textwidth,height=0.65\textwidth]{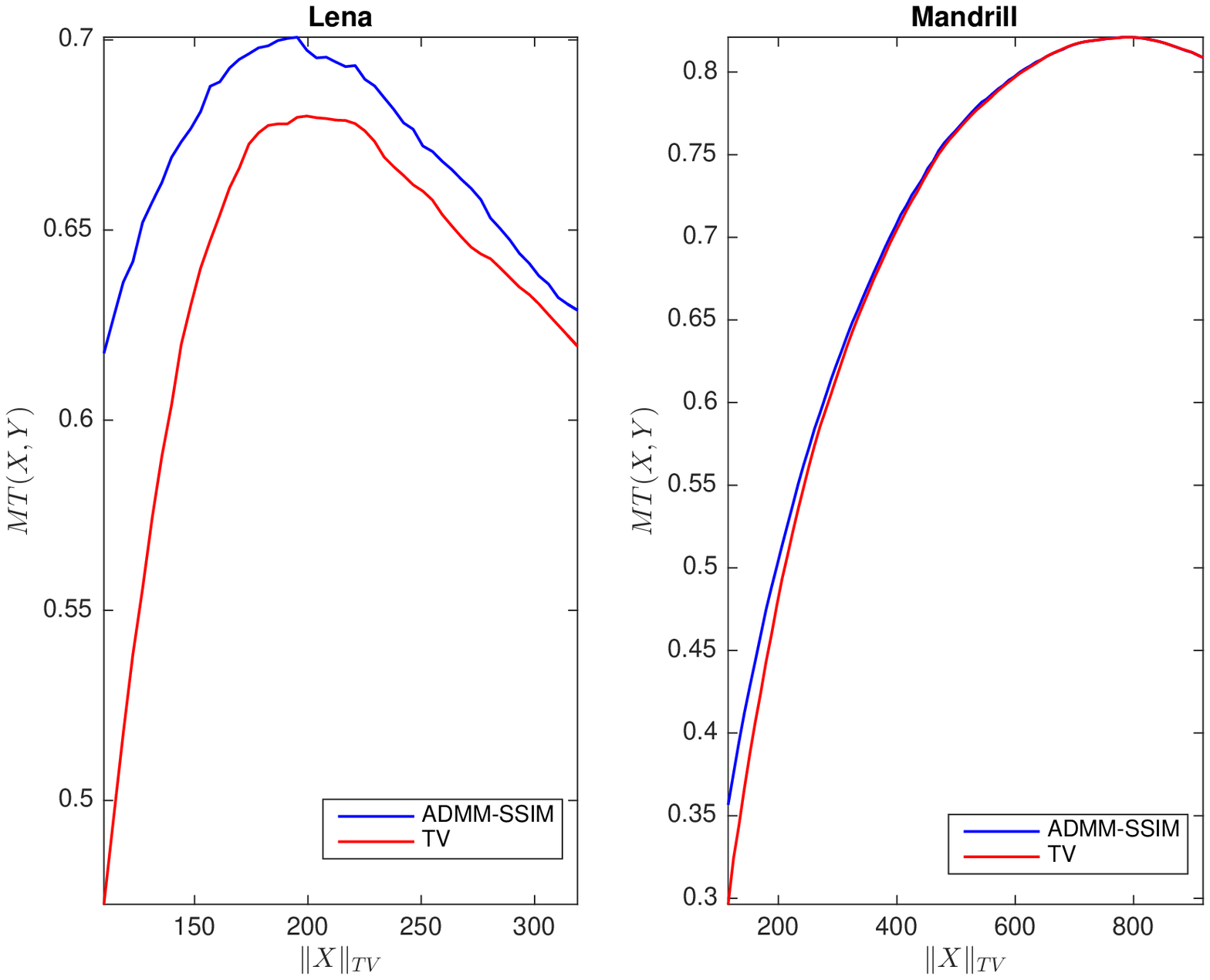}
\caption{The behaviour of the average SSIM of reconstructed images
	obtained from the proposed SSIM-based method and the classical
	$\ell_2$ method as a function of the TV seminorm of the reconstruction.
	{\bf Left:}  The \textit{Lena} image.  {\bf Right:} The
	\textit{Mandrill} image.
In the case of the \textit{Lena} image, the SSIM-based approach clearly outperforms the classical $\ell_2$ method.  For the \textit{Mandrill} image, however, the performance of both methods is, in general, very similar.}
\label{ssimtv}
\end{center}
\end{figure}

In order to have a general idea of the effectiveness of the SSIM-based methodology when regularization varies, in Figure \ref{ssimtv}, we show the behaviour of the MSSIM as a function of the TV seminorm of the reconstructions obtained by both the ADMM-SSIM and the TV approaches. The plot on the left shows the behaviour of the MSSIM for a noisy image patch of \textit{Lena} whereas the plot on the right shows the results for a corrupted image patch of \textit{Mandrill}. As expected, the plot on the right hand side shows that for images with low regularity---such as \textit{Mandrill}---the ADMM-SSIM and TV methods exhibit similar effectiveness over a wide range of regularization values.  On the other hand, for the image \textit{Lena}, one observes a significant difference between the performances of the two methods. This suggests that when strong regularization is required,  it is more advantageous to employ SSIM-based techniques over $\ell_2$ methods if certain visual features need to be recovered, provided that the reconstruction possesses some degree of regularity.

\subsubsection{Zooming}

Given a low resolution image $Y$, we wish to find an image $X$ which is a high resolution version of $Y$. This imaging task may be achieved by solving the following optimization problem \cite{ComboChambolle},
\begin{equation}
\min_{X} \, \|S(X)-Y\|_2^2+\lambda\|X\|_{TV},
\label{stv}
\end{equation}
where $S(\cdot)$ is a linear operator that consists of a blurring kernel followed by a subsampling procedure. Problem \eqref{stv} can be solved by different approaches such as the Chambolle-Pock algorithm \cite{ChambollePock} and ADMM \cite{BoydADMM}. In our particular case, we employed the ADMM method.

\noindent
The SSIM counterpart is given by
\begin{equation}
\min_{X} \, MT(S(X),Y)+\lambda\|X\|_{TV}.
\label{ztv}
\end{equation}
As mentioned in Section \ref{latv}, under certain circumstances, both Algorithm III and Algorithm IV can be used to perform the minimization of \eqref{ztv}.  In this study, we employed Algorithm IV.

\begin{figure}[htbp]
\begin{center}
\includegraphics[width=1\textwidth,height=0.7\textwidth]{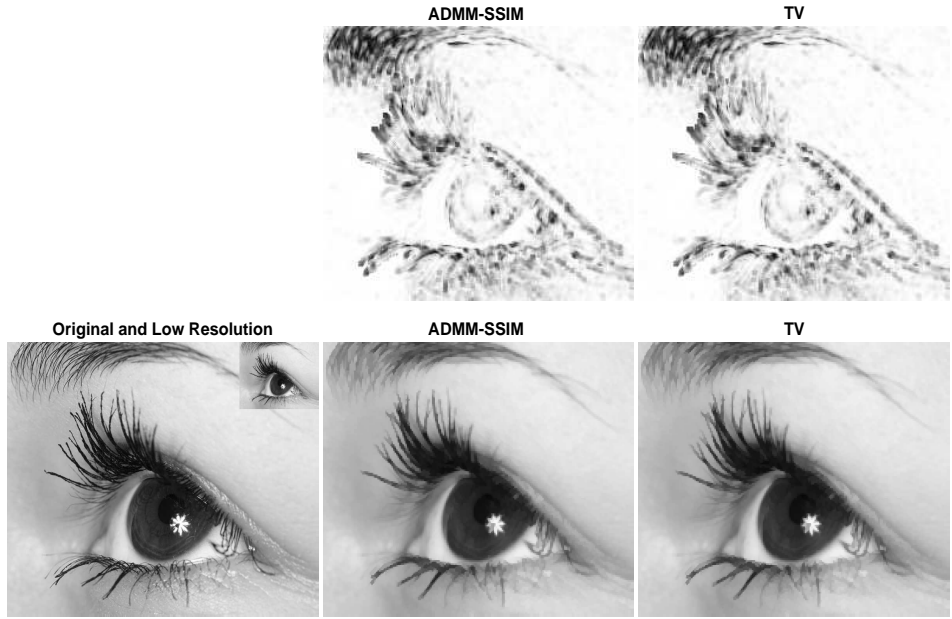}
\caption{Visual results for the Zooming application. In the first row, the SSIM maps are shown. Original image, low resolution observation, and reconstructions are presented in the bottom row. In this case, the zoom factor is four. The TV and the MSSIM for both reconstructions are 4563.26 and 0.7695 respectively. The image can be downloaded from \url{https://pixabay.com/fr/belle-gros-plan-oeil-sourcils-cils-2315/}.}
\label{zoomtv}
\end{center}
\end{figure}

Visual results are presented in Figure \ref{zoomtv}. In the first row, the SSIM maps are presented. The reconstructions yielded by both methods are presented in the bottom row along with the original image and its given low resolution counterpart. For the sake of a fair comparison, the strength of the regularization was adjusted so that the obtained reconstructions have the same TV. Also, in these experiments, the zoom factor is four.

\begin{figure}[htbp]
\begin{center}
\includegraphics[width=1\textwidth,height=0.7\textwidth]{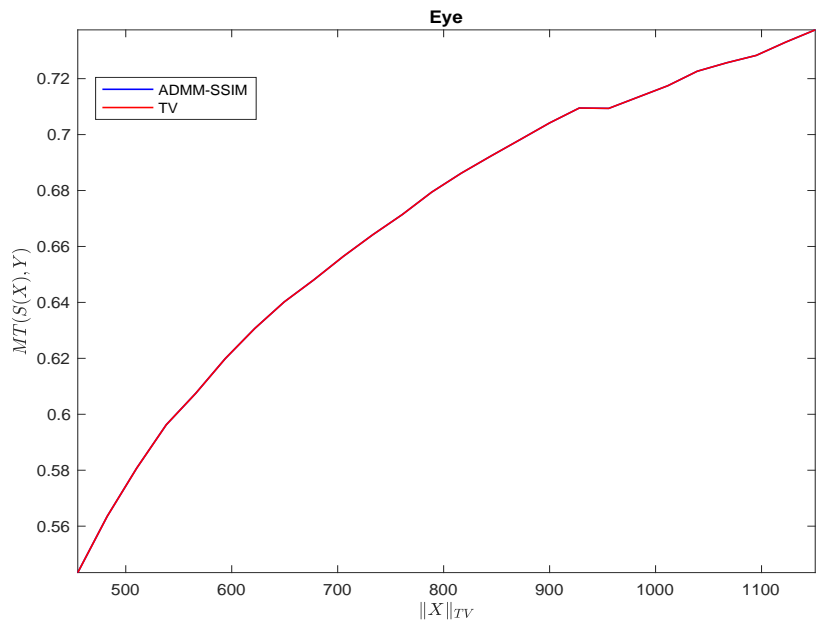}
\caption{A plot of the average SSIM as a function of the TV seminorm of the reconstructions yielded by both the ADMM-SSIM and $\ell_2$ methods. In these experiments, a patch of the reference image was employed for the zooming application. The performance of both methods is virtually the same.}
\label{tvzoom}
\end{center}
\end{figure}

It can be observed that both the proposed ADMM-SSIM method and the classical $\ell_2$ approach exhibit very similar performance. In fact, the MSSIM values for both reconstructions is 0.7695. In general, we found that the capabilities of both methods are virtually the same even when the strength of the regularization changes. This behaviour can be observed in Figure \ref{tvzoom}.

A possible explanation for this phenomenon is that many important features such as textures and fine details are either lost or barely present in the low resolution observation. The ADMM-SSIM method will obtain reconstructions which, after subsampling, will be visually as similar as possible to the low resolution observations.   These, in general, will not possess visual elements that are better captured by the SSIM-based approach. In other words, if, hypothetically speaking, we had two perfect reconstructions, one that is considered optimal with respect to the SSIM, and the other optimal with respect to $\ell_2$, both images would look almost identical when observed from afar.  This is equivalent to having the same low resolution image as an observation.

\subsubsection{Deblurring}

We perform the deblurring of an image by solving the following SSIM-based optimization problem:
\begin{equation}
\min_{X} \, MT(B(X),Y)+\lambda\|X\|_{TV},
\label{vtv}
\end{equation}
where $Y$ is the given blurred image, and $B(\cdot)$ is a blurring kernel. In all the experiments presented in this section, a Gaussian kernel was employed. Once again, we solve \eqref{vtv} by using Algorithm IV.  The $\ell_2$ counterpart of \eqref{vtv} is given by
\begin{equation}
\min_{X} \, \|B(X)-Y\|_2^2+\lambda\|X\|_{TV}.
\label{btv}
\end{equation}
The minimization of \eqref{btv} was performed using ADMM.

\begin{figure}[htbp]
\begin{center}
\includegraphics[width=1\textwidth,height=0.7\textwidth]{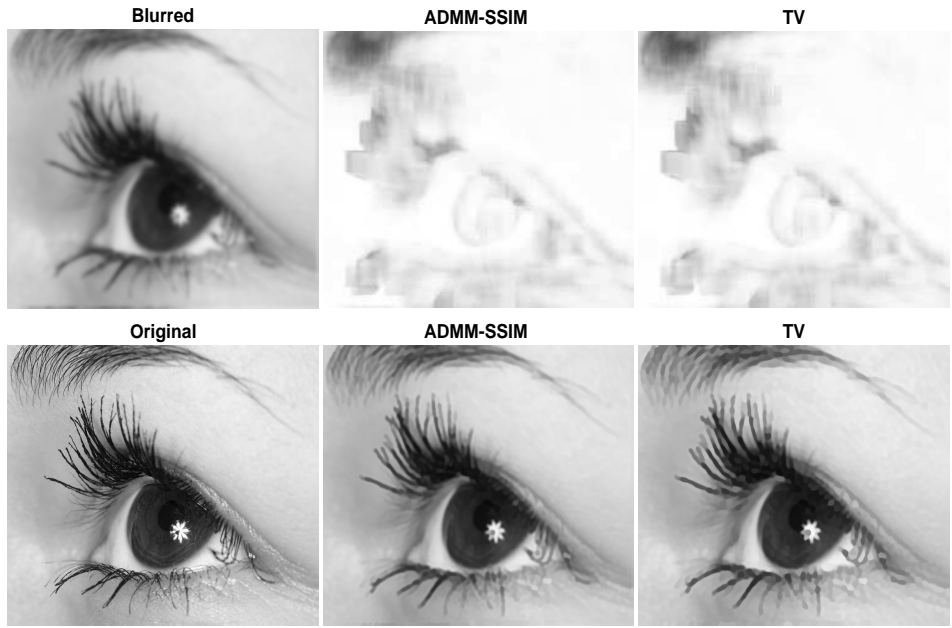}
\caption{Some visual results of the deblurring application. In the first row, the SSIM maps and the blurred observation are shown. Original image and reconstructions are presented in the bottom row. The standard deviation of the Gaussian blurring kernel was five. The MSSIM values of the reconstructions are 0.7517 and 0.7455 for the ADMM-SSIM and $\ell_2$ approaches, respectively. Notice that the fine features of the original image are better reconstructed by the proposed approach. The image can be downloaded from \url{https://pixabay.com/fr/belle-gros-plan-oeil-sourcils-cils-2315/}.}
\label{blurtv}
\end{center}
\end{figure}

Some visual results are shown in Figure \ref{blurtv}. The blurred observation along with the SSIM maps of the reconstructions are presented in the first row. The original image and the corresponding reconstructions of the two methods are presented in the bottom row. The TV of both reconstructions is 4000.  The standard deviation of the Gaussian blurring kernel was 5.

Overall, the effectiveness of both the $\ell_2$ and the ADMM-SSIM approaches is similar, however, the proposed method exhibits a better reconstruction of the fine details of the original image, e.g., the eyebrow, eyelashes and the reflection in the eye.  This can also be observed in the SSIM maps---compare, for example, the regions that correspond to the eyelashes of the upper eyelid. In this area, the SSIM map of the ADMM-SSIM reconstruction is brighter than the SSIM map of its $\ell_2$ counterpart. As for numerical results, the MSSIM values of the ADMM-SSIM and $\ell_2$ reconstructions are 0.7517 and 0.7455, respectively.

\section{Final Remarks}

In this paper, we have provided a general mathematical
as well as computational framework for
constrained and unconstrained SSIM-based optimization.
This framework provides a means of including SSIM as a fidelity term
in a wide range of applications.  Several
algorithms which can be used to accomplish a variety of
SSIM-based imaging tasks have also been proposed.
Problems in which both the $\ell_1$-norm and the
TV seminorm are used, either as constraints or as
regularization terms, have also been examined in this paper.  To the
best of our knowledge, these problems have been examined only to a small
degree in the literature \cite{Otero,Otero2014,YuShao,Brunet2017}.

We have employed a
simplified version of the SSIM which yields
a dissimilarity measure $T(x,y)$ that is, in fact, an example of a squared
normalized metric.
Mathematically, it is more desirable to work with $T(x,y)$ because of
its quasiconvex property.
Our SSIM-based optimization problems
are then formulated as minimization problems involving $T(x,y)$.

Formally, the dissimilarity measure $T(x,y)$ used in this paper
was obtained with the assumption that the vectors $x$ and $y$
have zero mean.  That being said,
some experimental results of Section \ref{ssimtvexp} suggest that
$T(x,y)$ may be used effectively in cases where the
the vectors have nonzero mean.

The results presented in this paper indicate that, in general, the performance of SSIM-based optimization schemes can be at least as good as that of their classical $\ell_2$ counterparts.  In some cases, SSIM-based methods appear to work better than $\ell_2$ approaches.  One example is the compression of images with low regularity (see Figure \ref{lemanssim}).  It appears that in these cases, the given (degraded) observation possesses at least partial information of the key features of the uncorrupted image.
That being said, the determination of problems for which
SSIM-based optimization outperforms
$\ell_2$-based methods, and {\em vice versa},
requires much more investigation.
The primary purpose of this paper, however, was to provide the framework
for such work.\\

\noindent
{\bf Acknowledgements:}  This work has been supported in part
by Discovery Grants (ERV and OM) from the Natural Sciences and
Engineering Research Council of Canada (NSERC).  Financial
support from the Faculty of Mathematics and the Department of
Applied Mathematics, University of Waterloo (DO) is also gratefully
acknowledged.

\section{Appendix}

\textbf{Proof of Theorem 4:} \textit{Let $f:X\subset\mathbb{R}^n\to\mathbb{R}$ be defined as in Eq. (\ref{thef}). Then, its Jacobian is Lipschitz continuous on any open convex set $\Omega\subset X$; that is, there exists a constant $L>0$ such that for any $x,w\in\Omega$,
\begin{equation}
\|J_f(x)-J_f(w)\|_F\leq L\|x-w\|_2 \, .
\end{equation}
Here, $\|\cdot\|_F$ denotes the Frobenius norm, and
\begin{equation}
L = C_1\|\Phi^T\Phi\|_F+\lambda (C_2\|z\|_2+C_3),~C_1,C_2,C_3>0.
\end{equation}}
\begin{proof}
Without loss of generality, and for the sake of simplicity, let the stability constant $C$ of the dissimilarity measure be zero. Also, let $y$ be a non-zero vector in $\mathbb{R}^m$. Let us define
\begin{equation}
s(x):=\frac{2x^Ty}{\|Dx\|_2^2+\|y\|_2^2},
\end{equation}
and
\begin{equation}
r(x):=\|Dx\|_2^2+\|y\|_2^2.
\end{equation}
Therefore, we have that $\|J_f(x)-J_f(w)\|_F$ is bounded by
\begin{eqnarray}
\|J_f(x)-J_f(w)\|_F&\leq&\|\Phi^T\Phi\|_F\|x\nabla s(x)^T-w\nabla s(w)^T\|_F+\nonumber\\
&&\lambda\|z(\nabla r(x)^T-\nabla r(w)^T)\|_F+\nonumber\\
&&\lambda\|x\nabla r(x)^T-w\nabla r(w)^T\|_F+\nonumber\\
&&|s(x)-s(w)|\|\Phi^T\Phi\|_F+\lambda|r(x)-r(w)|,
\end{eqnarray}
To show that $J_f$ is Lipschitz continuous on $\Omega$, we have to show that each term is Lipschitz continuous on $\Omega$ as well. Let us begin with the term $|r(x)-r(w)|$. By using the mean-value theorem for real-valued functions of several variables, we have that
\begin{equation}
|r(x)-r(w)|\leq\|2\Phi^T\Phi(\alpha x+(1-\alpha)w)\|_2\|x-w\|_2\\,
\end{equation}
for some $\alpha\in[0,1]$ and all $x,w\in\Omega$. Thus,
\begin{eqnarray}
|r(x)-r(w)|&\leq&2\|\Phi^T\Phi\|_2(\alpha\|x-w\|_2+\|w\|_2)\|x-w\|_2\\
&\leq&2\|\Phi^T\Phi\|_2(\|x-w\|_2+\|w\|_2)\|x-w\|_2.
\end{eqnarray}
Let $\sigma(\Omega)$ be the diameter of the set $\Omega$, that is,
\begin{equation}
\sigma(\Omega)=\sup_{x,v\in\Omega}\|x-v\|_2.
\end{equation}
Also, let $\rho(\Omega)$ be the $\ell_2$ norm of the largest element of the set $\Omega$, i.e.,
\begin{equation}
\rho(\Omega)=\sup_{x\in\Omega}\|x\|_2.
\end{equation}
Then,
\begin{eqnarray}
|r(x)-r(w)|&\leq&2\|\Phi^T\Phi\|_2(\sigma(\Omega)+\rho(\Omega))\|x-w\|_2\\
&\leq&K_1\|x-w\|_2,
\end{eqnarray}
where $K_1=2\|\Phi^T\Phi\|_2(\sigma(\Omega)+\rho(\Omega))$.

\noindent
As for $|s(x)-s(w)|$, in a similar fashion, we obtain that
\begin{equation}
|s(x)-s(w)|\leq\|\nabla s(\alpha x+(1-\alpha)w)\|_2\|x-w\|_2
\end{equation}
In fact, it can be shown that for any vector $v\in\mathbb{R}^n$, the norm of the gradient of $s$ is bounded by
\begin{equation}
\|\nabla s(v)\|\leq(\sqrt{2}+1)\frac{\|D\|_2}{\|y\|_2}.
\end{equation}
Let $K_2=(\sqrt{2}+1)\frac{\|D\|_2}{\|y\|_2}$. Thus, $|s(x)-s(w)|\leq K_2\|x-w\|_2$.

Regarding the term $\|x\nabla s(x)^T-w\nabla s(w)^T\|_F$, we have that the $ij$-th each entry of the $n\times n$ matrix $x\nabla s(x)^T-w\nabla s(w)^T$ is given by
\begin{equation}
\nabla_js(x)x_i-\nabla_js(w)w_i,
\end{equation}
where $\nabla_js(\cdot)$ is the $j$-th component of the gradient of $s(\cdot)$. By employing the mean-value theorem for functions of one variable we obtain that
\begin{equation}
|\nabla_js(x)x_i-\nabla_js(w)w_i|=\left|\frac{\partial}{\partial x_i}(\nabla_js(x(v)))\right||x_i-w_i|,
\end{equation}
for some $v\in\mathbb{R}$. Here, $x(v)=[x_1,\dots,x_{i-1},v,\dots,x_n]$. The partial derivative of the previous equation is bounded, which can be proved using the classical triangle inequality and differential calculus. Given this, we have that
\begin{eqnarray}
\left|\frac{\partial}{\partial x_i}(\nabla_js(x))(v)\right|&\leq&(\sqrt{2}+3)\frac{\|\Phi_i^T\|_2\|\Phi_j^T\|_2}{\|y\|_2^2}+(2\sqrt{3}+2)\frac{\|\Phi_j^T\|_2}{\|y\|_2^3}\\
&=&K_{ij},
\end{eqnarray}
where $\Phi_k^T$ is the $k$-th row of the the transpose of the matrix $\Phi$. Therefore,
\begin{equation}
|\nabla_js(x)x_i-\nabla_js(w)w_i,|\leq K_{ij}|x_i-w_i|.
\end{equation}
Using this result, we can conclude that
\begin{equation}
\|x\nabla s(x)^T-w\nabla s(w)^T\|_F\leq K_3\|x-w\|_2,
\end{equation}
where $K_3$ is equal to
\begin{equation}
K_3=n\max_{1\leq i,j\leq n}K_{ij};
\end{equation}
that is, $K_3$ is equal to the largest $K_{ij}$ times $n$.

\noindent
In a similar manner, it can be shown that
\begin{equation}
\|x\nabla r(x)^T-w\nabla r(w)^T\|_F\leq K_4\|x-w\|_2,
\end{equation}
where $K_4$ is given by
\begin{equation}
K_4=\max_{1\leq i,j\leq n}\{2nK_1\|\Phi_j^T\|_2(\|\Phi_i^T\|_2+\|\Phi\|_2)\}.
\end{equation}

\noindent
As for the term $\lambda\|z(\nabla r(x)^T-\nabla r(w)^T)\|_F$, this is equal to
\begin{equation}
\lambda\|z(\nabla r(x)^T-\nabla r(w)^T)\|_F=2\lambda\|z(\Phi^T\Phi(x-w))^T\|_F.
\end{equation}
Each $ij$-th entry of the matrix $z(\Phi^T\Phi(x-w))^T$ is given by $z_i(\Phi^T\Phi_j(x-w))^T$, where $\Phi^T\Phi_j$ is the $j$-th row of $\Phi^T\Phi$. Then, we have that
\begin{equation}
2\lambda\|z(\Phi^T\Phi(x-w))^T\|_F=2\lambda\sqrt{\sum_i^n\sum_j^n|z_i(\Phi^T\Phi_j(x-w))^T|^2}.
\end{equation}
Therefore,
\begin{eqnarray*}
2\lambda\sqrt{\sum_i^n\sum_j^n|z_i(\Phi^T\Phi_j(x-w))^T|^2}&\leq&2\lambda\sqrt{\sum_i^n\sum_j^n|z_i|^2\|\Phi^T\Phi_j\|_2^2\|x-w\|_2^2}\\
&\leq&2\lambda\sqrt{\sum_i^n|z_i|^2\sum_j^n\|\Phi^T\Phi_j\|_2^2}\|x-w\|_2\\
&=&2\lambda\|z\|_2\sqrt{\sum_j^n\|\Phi^T\Phi_j\|_2^2}\|x-w\|_2\\
&=&\lambda K_5\|z\|_2\|x-w\|_2\\,
\end{eqnarray*}
where $K_5=2\sqrt{\sum_j^n\|\Phi^T\Phi_j\|_2^2}$.

\noindent
Finally, we obtain that
\begin{equation*}
\|J_f(x)-J_f(w)\|_F\leq[(K_2+K_3)\|\Phi^T\Phi\|_F+\lambda(K_5\|z\|_2+K_1+K_4)]\|x-z\|_2,
\end{equation*}
which completes the proof.
\end{proof}

\end{document}